\documentclass{elsarticle}

\usepackage{amsmath}
\usepackage{amssymb}
\usepackage{bm}
\usepackage{tikz}
\usepackage{pgfplots}
\usepgflibrary{arrows}
\usepackage{rotating}
\usepackage{comment,xspace}
\graphicspath{{figs/}}

\usepackage{hyperref}

\newtheorem{thm}{Theorem}

\newproof{pf}{Proof}

\journal{Journal of Computational and Applied Mathematics}

\begin{document}

\begin{frontmatter}

\title{High order numerical schemes for solving  fractional powers of elliptic operators}

\author[vgtu]{Raimondas {\v C}iegis}
\ead{raimondas.ciegis@vgtu.lt}

\author[nsi,nefu]{Petr N. Vabishchevich\corref{cor}}
\ead{vabishchevich@gmail.com}

\address[vgtu]{Vilnius Gediminas Technical University, Sauletekio al. 11, 
10223 Vilnius,  Lithuania}
\address[nsi]{Nuclear Safety Institute, Russian Academy of Sciences, 52, B. Tulskaya, Moscow, Russia}
\address[nefu]{North-Eastern Federal University, 58, Belinskogo, Yakutsk, Russia}

\cortext[cor]{Corresponding author}

\begin{abstract}
In many recent applications when new materials and technologies are developed it is important 
to describe and simulate  new nonlinear and nonlocal diffusion transport processes. A general 
class of such models deals with nonlocal fractional power elliptic operators. In order 
to solve these problems numerically it is proposed    
(Petr N. Vabishchevich, Journal of Computational Physics. 2015, Vol. 282, No.1, pp.
289--302) to consider equivalent local nonstationary initial value  pseudo-parabolic 
problems.  Previously such problems were solved by using the standard implicit  backward 
and  symmetrical Euler  methods. In this paper we use the one-parameter family of
three-level finite difference 
schemes for solving the initial value problem for the  first order nonstationary 
pseudo-parabolic problem. The fourth-order approximation scheme is developed by
selecting the optimal value of the weight parameter. The results of the theoretical 
analysis are supplemented by results of extensive computational experiments.        
\end{abstract}

\begin{keyword}
elliptic operator \sep fractional power of an operator \sep finite element approximation \sep
three-level schemes \sep stability of difference schemes

\MSC[2010] 26A33 \sep 35R11 \sep 65F60 \sep 65M06
\end{keyword}

\end{frontmatter}

\section{Introduction}
In many recent applications the new mathematical models are proposed, which are based 
on fractional derivative equations in time and space coordinates 
\cite{baleanu2012fractional,eringen2002nonlocal,kilbas2006theory}. Very different 
applied mathematical models of physics, biology or finance describe a  subdiffusion 
(represented by fractional-in-time derivatives) or superdiffusion (represented by 
fractional-in-space derivatives) models. The latter problems are often simulated by
using fractional power
elliptic operators. 

Different numerical techniques, such as finite difference, finite volume methods, can 
be used  to approximate problems with  fractional power elliptic operators. In this
paper we will use 
the method of finite elements, since this method is well-suited  to solve problems 
in non-regular domains and to use non-uniform  adaptive grids
\cite{KnabnerAngermann2003,QuarteroniValli1994}. The implementation of such 
algorithms require to   
compute
the action of a matrix (operator)
function on a vector
$\varPhi(A) b$, where $A$ is a given 
matrix (operator) and $b$ is a given vector. For example, in order to compute 
the solution of the discrete fractional order  elliptic problem, we get  
$\varPhi(z) = z^{-\alpha}$, where  $0 < \alpha < 1$.
There exist various approaches how to compute
$\varPhi(A) b$
\cite{higham2008functions}.

The most important class of iterative methods for this purpose are
Krylov subspace methods. They are used to solve systems of linear equations
obtained after approximation of fractional power elliptic problems (see, e.g. 
\cite{ilic2009numerical}).
A comparison of different approaches to solve fractional-in-space reaction-diffusion
equations is done   
in \cite{burrage2012efficient}. In particular the integral and  
adaptively preconditioned Lanczos method are analyzed.

The most straightforward algorithm to solve such systems  is to construct explicitly
eigenvectors  and eigenvalues of the given discrete elliptic operator and to diagonalize 
the matrix $A$
\cite{ilic2005numerical,ilic2006numerical,bueno2012fourier}.
But we  should note that
the direct implementation of this approach is very expensive for general elliptic
operators in multidimensional domains. It requires the
computation of all eigenvectors and eigenvalues of very large matrices.

A general approach to solve fractional power elliptic problems is based on some approximation  
of the nonlocal operator.  

One can adopt a general approach to solve numerically equations involving fractional power
of operators by
a popular method is to split  the task to solve numerically equations involving 
fractional power into two steps.
First the original elliptic operator is approximated and then the fractional power
of its discrete variant is taken. 
Using Dunford-Cauchy formula the elliptic operator is represented
as a contour integral in the complex plane. Then applying appropriate quadratures with
integration nodes in the complex plane we get a method that involves only
inversion of the original elliptic operator. 
The approximate operator is treated as a sum of resolvents
\cite{gavrilyuk2004data,gavrilyuk2005data},
ensuring the exponential convergence of quadrature approximations.

In paper
\cite{bonito2015numerical} a more promising quadratures algorithm is proposed, when the 
integration nodes are selected in the real axis. The new method is based on the 
integral representation of the power operator  
\cite{krasnoselskii1976integral}. In this case the inverse operator of the fractional 
power elliptic problem 
is treated as a sum of inverse operators of elliptic operators. 

Such a rational approximation is obtained when
the fractional power of the operator is 
approximated  by using the 
Gauss-Jacobi quadrature formulas for the corresponding integral representation.
In this case, we have (see \cite{frommer2014efficient,AcetoNovat}) a Pade-type 
approximation of the power function with
a fractional exponent. The optimal rational approximations are investigated in
\cite{vciegis2017parallel,harizanov2018optimal}.

A separate class of methods approximates the solution of fractional power elliptic
problem by some auxiliary problem of high dimension. In 
\cite{Caffarelli} it is shown that the solution of the
fractional {L}aplacian problem can be obtained as a solution of the elliptic 
problem on the semi-infinite cylinder domain. This idea is used to construct numerical 
algorithms for solving stationary and non-stationary problems with fractional power 
elliptic operators \cite{nochetto2015pde,nochetto2016pde}.

In \cite{vabishchevich2014numerical}, for solving 
fractional power elliptic problems
we have proposed a numerical algorithm 
on the basis of a transition to a pseudo-parabolic equation, so called 
Cauchy problem method.
The computational algorithm is simple for practical use, robust, 
and applicable to solving
a wide class of problems. We have used this algorithm also for solving
the nonstationary problem with fractional power elliptic operators 
\cite{vabishchevich2016cnumerical}.

For the auxiliary Cauchy problem, standard two-level schemes are applied. 
Depending on the weight parameters the first and second order 
accuracy of the approximation is obtained.  For many 
applied problems a small number of pseudo-time steps is sufficient to
get a good approximation of the solution of the discrete fractional equation.
The efficiency of this algorithm is improved in \cite{duan2018numerical},
where a special graded grid in pseudo-time is used.    

Another possibility to increase the accuracy of approximations is to use high order 
discrete schemes for solving the auxiliary pseudo-parabolic equation.  
In this paper we propose and  investigate a fourth order three-level scheme.  

The paper is organized as follows. In Section 2 
a problem for a fractional power
of elliptic operator is formulated.
In Section 3 the Cauchy problem method is given. The main results are described in Section 4,
where unconditionally stable fourth-order  three-level scheme is proposed and investigated. 
Section 4 provides results of computational experiments, they illustrate the theoretical 
results on the approximation accuracy  of fractional power problems. A model two dimensional
problem is solved by using different numerical schemes.  
At the end of the work the main results of our study are summarized.

\section{Problem Formulation}

In a bounded domain
$\Omega \subset {\mathbb R}^d$, $d=2,3$ 
with the Lipschitz continuous boundary
$\partial\Omega$
we solve the boundary value problem for the fractional power 
elliptic operator.
The following elliptic  operator is defined by:
\begin{equation}\label{1}
  \mathcal{A} u = - {\rm div}  ( k({\bm x}) {\rm grad} u ) + c({\bm x}) u
\end{equation}
where
$0 < k_1 \leq k({\bm x}) \leq k_2$, $c({\bm x}) \geq c_0 > 0$.
On $\partial\Omega$ 
the functions $u({\bm x})$ satisfy the boundary conditions 
\begin{equation}\label{2}
  k({\bm x}) \frac{\partial u }{\partial n } + \mu ({\bm x}) u = 0,
  \quad {\bm x} \in \partial \Omega ,
\end{equation}
where  $\mu ({\bm x}) \geq 0, \  {\bm x} \in \partial \Omega$.

In the Hilbert space
$H = L_2(\Omega)$ we define the scalar product and norm:
\[
  (u,v) = \int_{\Omega} u({\bm x}) v({\bm x}) d{\bm x},
  \quad \|u\| = (u,u)^{1/2} .
\]
Next we introduce the eigenvalue problem \cite{bookEvans} for  
\eqref{1}, \eqref{2}:
find $\varphi_j \in H$ and $\lambda_j \in \mathbb{R}$ so that
\[
 \mathcal{A}  \varphi_k = \lambda_k \varphi_k, 
 \quad \bm x \in \Omega , 
\]
\[
  k({\bm x}) \frac{\partial  \varphi_k}{\partial n } + \mu ({\bm x}) \varphi_k = 0,
  \quad {\bm x} \in \partial \Omega.  
\]
The eigenvectors are numbered in such a way, that
\[
 \lambda_1 \leq \lambda_2 \leq ... .
\]
This spectral problem has full set of eigenfunctions 
$
\{ \varphi_k, \;  \|\varphi_k\| = 1, \;  k = 1,2, ...  \} 
$
that span the space
$L_2(\Omega)$:
\[
 u = \sum_{k=1}^{\infty} (u,\varphi_k) \varphi_k .
\]
We assume that the operator
 $\mathcal{A}$ is defined on the domain
\[
 D(\mathcal{A} ) = \big\{ u:  \; u(\bm x) \in L_2(\Omega), \;\;
 \sum_{k=0}^{\infty} | (u,\varphi_k) |^2 \lambda_k < \infty \big\} .
\]
Then 
$\mathcal{A}$
is self-adjoint and coercive
\begin{equation}\label{3}
  \mathcal{A}  = \mathcal{A}^* > \delta \mathcal{I} ,
  \quad \delta > 0 ,    
\end{equation}
where $\mathcal{I}$ is the identity operator in $H$. For  $\delta$ 
we have $\delta < \lambda_1$.
In most applications, the value of 
$\lambda_1$ is unknown and it should be obtained numerically by solving 
the eigenvalue problem. In our analysis we assume that a reliable 
positive bound from below is known
 $\delta < \lambda_1$ in  \eqref{3}. 

The fractional power of
$\mathcal{A}$ is defined by 
\[
 \mathcal{A}^{\alpha} u =  \sum_{k=0}^{\infty} (u,\varphi_k) \lambda_k^{\alpha }  \varphi_k ,
\]
where $0 < \alpha < 1$.
Now we define the boundary value problem for the fractional power of
$\mathcal{A}$. The solution
$u(\bm x)$ satisfies the equation
\begin{equation}\label{4}
  \mathcal{A}^{\alpha} u = f .
\end{equation}

We approximate the problem
\eqref{4} by using the finite element method
\cite{brenner2008mathematical}.
For the elliptic problem \eqref{1}, \eqref{2}
the bilinear form is defined by
\[
 a(u,v) = \int_{\Omega } \left ( k \, {\rm grad} \, u \, {\rm grad} \, v + c \, u v \right )  d {\bm x} +
 \int_{\partial \Omega } \mu \, u v d {\bm x} .
\]
Due to \eqref{3}, we have that
\[
a(u,u) > \delta \|u\|^2 .  
\]
We consider a standard sub-space of finite elements
$V_h \subset H^1(\Omega)$. Let us consider a triangulation of the domain 
$\Omega$ into triangles and let
$\bm x_i$,  $i = 1,2, ..., M_h$ be  the vertexes of these triangles. 
As a nodal basis we take the functions
$\chi_i(\bm x) \subset V_h$,  $i = 1,2, ..., M_h$:
\[
 \chi_i(\bm x_j) = \left \{
 \begin{array}{ll} 
 1, & \mathrm{if~}  i = j, \\
 0, & \mathrm{if~}  i \neq  j .
 \end{array}
 \right . 
\]
Then for $v \in V_h$ we have
\[
 v(\bm x) = \sum_{i=i}^{M_h} v_i \chi_i(\bm x),
\]
where $v_i = v(\bm x_i), \ i = 1,2, ..., M_h$.
We define  the discrete elliptic operator $A$
\[
 a(u,v) = \ (Au, v),
 \quad \forall \ u,v \subset V_h.  
\]
Similar to 
\eqref{3}, the following estimates are valid for $A$:
\begin{equation}\label{5}
 A = A^* > \delta I ,
 \quad \delta > 0 .  
\end{equation}
The corresponding finite element approximation of equation 
\eqref{4} is: find
$v \subset V_h$
\begin{equation}\label{6}
 A^\alpha v = \psi , 
\end{equation}
where $\psi = P f$ and 
$P$ is the
$L_2$ projection on 
$V_h$.
In view of \eqref{5}, for the solution \eqref{6} we get
the following simple a priori estimate:
\begin{equation}\label{7}
 \|v\| \leq \delta^{-\alpha} \|\psi\|.
\end{equation}

\section{Cauchy problem method}

For solving numerically problem \eqref{6} we use the Cauchy problem method, 
proposed in
\cite{vabishchevich2014numerical}. This method is based on 
the equivalence of \eqref{6} to  
an auxiliary pseudo-time
evolutionary problem.
Assume that
\[
 w(t) = \delta^{\alpha } (t (A - \delta I) + \delta I)^{-\alpha} w(0) .
\]
Therefore
\[
 w(1) =  \delta^{\alpha} A^{-\alpha} w(0)
\]
and then $v = w(1)$ if $w(0) = \delta^{-\alpha} \psi$.
The function $w(t)$ satisfies the evolutionary equation
\begin{equation}\label{8}
  B(t) \frac{d w}{d t} + D w = 0 ,
  \quad 0 < t \leq 1 ,
\end{equation}
where
\[
 B = \frac{1}{\alpha} (t D + \delta I) ,
 \quad D = A - \delta I .
\]
We supplement \eqref{8} with the initial condition
\begin{equation}\label{9}
 w(0) = \delta^{-\alpha} \psi.  
\end{equation}
By \eqref{5}, we get
\begin{equation}\label{10}
 D = D^* > 0 .
\end{equation}
The solution of equation \eqref{6} can be defined as the solution of the Cauchy problem
\eqref{8}, \eqref{9} at the final pseudo-time moment $t=1$.

For the solution of the problem \eqref{8}, \eqref{9}, it is possible to obtain various a priori estimates.
Here we restrict only to a simple
estimate that is consistent with the estimate \eqref{7}: 
\begin{equation}\label{11}
  \|w(t)\|  \leq \|w(0)\| .
\end{equation}
In order to prove 
\eqref{11}, it is sufficient to  multiply scalarly equation \eqref{8} 
by $\alpha w + t dw/ dt$.

To solve numerically the problem \eqref{8}, \eqref{9},
the simple implicit two-level Euler scheme can be used  \cite{Samarskii1989}.
Let $\tau$ be the step-size of a uniform grid in time: 
\[
w_n = w(t_n), \ t_n = n \tau, \quad n = 0,1, ..., N, \quad  N\tau = 1.
\]
Let us approximate equation \eqref{8} by the implicit two-level scheme
\begin{equation}\label{12}
(t_{n+\sigma} D + \delta I) \frac{ w_{n+1} - w_{n}}{\tau }
 + \alpha D w_{n+\sigma} = 0,  \quad n = 0,1, ..., N-1,
\end{equation} 
\begin{equation}\label{13}
 w_0 = \delta^{-\alpha} \psi .
\end{equation}
We use the notation
\[
  t_{n+\sigma} = \sigma t_{n+1} + (1-\sigma) t_{n},
  \quad w_{n+\sigma} = \sigma w_{n+1} + (1-\sigma) w_{n}.
\]
For sufficiently smooth $w(t)$ and
$\sigma =0.5$ (the Crank-Nicolson type scheme), the difference 
scheme \eqref{12}, \eqref{13} approximates the problem
\eqref{8}, \eqref{9}
with the second order, and with the first order for all 
 other values of $\sigma$.

\begin{thm}\label{t1}
For $\sigma \geq 0.5$ the difference scheme \eqref{12}, \eqref{13}
is unconditionally stable with respect to the initial data.
The approximate solution satisfies the estimate
\begin{equation}\label{14}
  \|w_{n+1}\|  \leq \delta^{-\alpha} \psi  , 
  \quad n = 0,1, ..., N-1.
\end{equation}
\end{thm}
\begin{pf}
In order to prove this result we 
rewrite equation \eqref{12} in the following form:
\[
 \delta \frac{ w_{n+1} - w_{n}}{\tau } + D \left( \alpha w_{n+\sigma} + t_{n+\sigma}\frac{ w_{n+1} - w_{n}}{\tau } \right ) = 0 . 
\]
Multiplying scalarly it by
\[
 \alpha w_{n+\sigma} + t_{n+\sigma}\frac{ w_{n+1} - w_{n}}{\tau },
\]
and due to \eqref{10}, we get
\[
 \left ( \frac{ w_{n+1} - w_{n}}{\tau }, w_{n+\sigma} \right )  \leq 0 . 
\]
Since
\[
 w_{n+\sigma} = \left( \sigma - \frac{1}{2} \right ) \tau \frac{w_{n+1} - w_{n}}{\tau }   +
 \frac{1}{2} (w_{n+1} + w_{n})
\]
then for $\sigma \geq 0.5$ we get the inequalities 
\[
  \|w_{n+1}\|  \leq \|w_n\|  , 
  \quad n = 0,1, ..., N-1 .  
\]
Applying these estimates recursively we prove the validity of \eqref{14}.
\end{pf}

\section{Three level schemes}

In this section we consider high order schemes. They are based on three
level finite difference schemes. For solving problem \eqref{8}, \eqref{9}
we use the symmetrical scheme 
\begin{equation}\label{15}
\begin{split}
 B_n \frac{w_{n+1} - w_{n-1}}{2 \tau } & + D (\sigma w_{n+1} + (1-2 \sigma) w_{n} + \sigma w_{n-1}) = 0, \\
 & \quad n = 1,2, \ldots, N-1 ,
\end{split}
\end{equation}
with the given initial conditions
\begin{equation}\label{16}
 w_0 = \delta^{-\alpha} \psi, 
 \quad  w_1 = \overline{w}_1 .
\end{equation}
We note that  $\overline{w}_1$ should be computed by applying some 
two level numerical algorithm and the accuracy of this approximation 
should be the same as of the main scheme \eqref{15}. More details
will be given below. 

It is well-known that for sufficiently smooth solutions
the symmetrical sche\-me \eqref{15} approximates problem 
\eqref{8}--\eqref{10} with the second order accuracy. 

Next we formulate the stability conditions for 
the scheme \eqref{15}, \eqref{16}. 
Here we use the general stability results for operator-difference  schemes
\cite{Samarskii1989,SamarskiiMatusVabischevich2002}. 

Let $S$ be a self-adjoint positive operator in $H$. Then we introduce
the new Hilbert space
$H_S$, generated by operator $S$, it consists of elements from $H$ 
equipped with the energy norm 
\[
 (y,v)_S = (Sy,v) ,
 \quad \| y \|_S = (Sy,y)^{1/2}.
\]

\begin{thm}\label{t-2}
For $\sigma > 0.25$ the three-level scheme \eqref{15}, \eqref{16}
is unconditionally stable with respect to the initial data.
The approximate solution satisfies the estimate
\begin{equation}\label{17}
  \mathcal{E}_{n+1} \leq  \mathcal{E}_{n} ,
  \quad n = 1,2, ..., N-1 ,
\end{equation}
where
\begin{equation}\label{18}
 \mathcal{E}_{n} = \frac{1}{4} \|w_{n} + w_{n-1}\|_D^2
 + \left (\sigma - \frac{1}{4} \right ) \|w_{n} - w_{n-1}\|_D^2 .
\end{equation}
\end{thm}
\begin{pf}
We rewrite equation \eqref{15} in the following form
\begin{equation}\label{19}
 B_n \frac{w_{n+1} - w_{n-1}}{2 \tau } + 
\sigma D (w_{n+1} - 2 w_{n} + w_{n-1}) + A w_{n}= 0.
\end{equation}
let us introduce new functions 
\[
 y_n  = \frac{1}{2} (w_{n} + w_{n-1}),
 \quad z_n = w_{n} - w_{n-1}.
\]
Taking into account that
\[
 w_{n} = \frac{1}{4} (w_{n+1} + 2 w_{n} + w_{n-1}) - \frac{1}{4} (w_{n+1} - 2 w_{n} + w_{n-1})
\]
we write 
\eqref{19} in the following form 
\begin{equation}\label{20}
 B_n \frac{z_{n+1} + z_{n}}{2 \tau } + \left (\sigma - \frac{1}{4} \right ) D (z_{n+1} - z_{n})
 + D \frac{y_{n+1} + y_{n}}{2}  = 0 .
\end{equation}
Multiplying it by
\[
 2 (y_{n+1} - y_{n}) = z_{n+1} + z_{n}
\]
and taking a discrete inner product,  in view
that $D$ is self-adjoint, we get
\[
\begin{split}
 \frac{1}{2 \tau} (B_n(z_{n+1} + z_{n}), z_{n+1} + z_{n}) & +
 \left (\sigma - \frac{1}{4} \right )  \big ( (D z_{n+1}, z_{n+1}) - (D z_{n}, z_{n}) \big ) \\
 & + (D y_{n+1}, y_{n+1}) -  (D y_{n}, y_{n}) = 0 .
\end{split}
\]
Since
$B_n > 0$,
then we easily get the estimates
\eqref{17}.
\end{pf}

Next we consider how to define the initial condition 
\eqref{16} for $w_1$. A general approach is  
to use some two-level 
solver for $t \in [0,\tau]$. For example, it is possible
to apply the symmetrical scheme \eqref{12}, with $\sigma = 0.5$: 
\begin{equation}\label{21}
 B_{1/2} \frac{w_{1} - w_{0}}{\tau } + D \frac{w_{1} + w_{0}}{2 } = 0 .
\end{equation}
It follows from Theorem~\ref{t1}, that the scheme \eqref{21} is 
unconditionally stable  
\begin{equation}  \label{22}
  \|w_1\| \leq \|w_0\|
\end{equation}
and its approximate solution $w_1$ converges to $w(t_1)$ with second order.

It is interesting to see if some explicit schemes can be used 
to find the initial condition for $w_1$. One possibility 
is to consider the explicit forward 
Euler scheme 
\begin{equation}\label{23}
 \frac{w_{1} - w_{0}}{\tau } + \frac{\alpha}{\delta} D w_{0} = 0, 
\end{equation}
here the equality $B(0) = \alpha^{-1} \delta I$ is used. In general 
for sufficiently smooth solutions the $O(\tau^2)$ accuracy is expected
for $w_1$. Let us denote the error of the solution of \eqref{23} 
$\tilde z_n = w_n - w(t_n)$, $n=0,1$. The function $z_n$ satisfies the  
equation   
\[
 \frac{\tilde z_{1} - \tilde z_{0}}{\tau } + \frac{\alpha}{\delta} 
D \tilde z_{0} = \psi_1, 
\]
where $\psi_1$ is the standard approximation error. Since $\tilde z_0 = 0$,
then with sufficiently small step $\tau$ the error $\tilde z_1$ can be estimated as 
\[
    \Vert \tilde z_1 \Vert  \leq \tau \Vert \psi_1 \Vert.
\] 
For a sufficiently smooth solution of \eqref{8}, we have that 
$\Vert \psi_1 \Vert \leq C \tau$.

More interesting second order explicit schemes can be constructed 
by using the well-known method
described in \cite{Samarskii1989}. The accuracy of the
basic  forward 
Euler scheme \eqref{23} is increased by using the differential properties
of the solution of equation \eqref{8}   
\begin{equation}\label{24}
 \frac{w_{1} - w_{0}}{\tau } + \frac{\alpha}{\delta} D w_{0} - \frac{\alpha(1+\alpha)}{2 \delta^2} \tau  D^2 w_{0} = 0 .
\end{equation}
We rewrite \eqref{24} in the following form
\begin{equation}\label{25}
 w_1 = R w_0,
 \quad R = I - \frac{\alpha}{\delta} \tau D 
+ \frac{\alpha(1+\alpha)}{2 \delta^2} \tau^2  D^2.
\end{equation}
Then the stability estimate \eqref{22} is valid if $\|R\| \leq 1$. 
For a self-adjoint operator $R$ this estimate is equivalent to
the following two-side
estimates
\begin{equation}\label{26}
 - I \leq R \leq I .
\end{equation}
Due to \eqref{25} the right inequality of \eqref{26} can be written as
\[
 \frac{1+\alpha}{2 \delta} \tau I
 \leq  D .
\]
Then the following restrictions on the time step are obtained  
\begin{equation}\label{27}
 \tau \leq \tau_0 = \frac{2 \delta}{1+\alpha} \frac{1}{\|D\|} .
\end{equation}
The left inequality  can be rewritten as
\[
 I - \frac{\alpha}{2\delta} \tau D + \frac{\alpha(1+\alpha)}{4 \delta^2} \tau^2 \geq 0.
\]
For the given values of powers $\alpha$ this inequality is always valid, since 
\[
 I - \frac{\alpha}{2\delta} \tau D + \frac{\alpha(1+\alpha)}{4 \delta^2} \tau^2 D^2
   = \left ( I  - \frac{\alpha}{4\delta} \tau D \right )^2 +
   \frac{\alpha(4+3\alpha)}{16 \delta^2} \tau^2 D^2 .
\]
Due to the obtained stability restrictions \eqref{27} the explicit scheme
\eqref{24} is not recommended for  solving real applications.   

It is important to note that in the family of second order unconditionally
stable three-level 
schemes \eqref{15}, \eqref{16} it is possible to find such a value 
of the parameter $\sigma$ which leads to the high order accuracy 
scheme.  

Using the Taylor expansions we get the relations
\[
 w_{n+1} + w_{n-1} = 2 w+\tau^2 \frac{d^2 w}{d t^2}+\mathcal{O} (\tau^4),
\]
\[
 \frac{w_{n+1} - w_{n-1}}{2\tau } = \frac{d w}{d t} + \frac{\tau^2}{6} \frac{d^3 w}{d t^3} + \mathcal{O} (\tau^4) .
\]
Then the residual of the scheme can be written as 
\[
\begin{split}
 B_n \frac{w_{n+1} - w_{n-1}}{2\tau } & + D (\sigma w_{n+1} + (1-2\sigma) w_n + \sigma w_{n-1}) \\
 & =  B \frac{w_{n+1} - w_{n-1}}{2\tau } + (1-2\sigma) D u_w + \sigma D (w_{n+1} +  w_{n-1}) \\
 & = B \frac{d w}{d t} +  \frac{\tau^2}{6} B \frac{d^3 w}{d t^3} + D w + \sigma \tau^2 D \frac{d^2 w}{d t^2}
 + \mathcal{O} (\tau^4) .
\end{split}
\]
Taking the solution of
\eqref{8} we get
\begin{equation}\label{28}
\begin{split}
 B_n \frac{w_{n+1} - w_{n-1}}{2\tau } & + D (\sigma w_{n+1} + (1-2\sigma) w_n + \sigma w_{n-1}) \\
 & =\frac{\tau^2}{6} B \frac{d^3 w}{d t^3} + \sigma \tau^2 D \frac{d^2 u}{d t^2}
 + \mathcal{O} (\tau^4) .
\end{split}
\end{equation}
Differentiation of the equation  \eqref{8} leads to the equality
\[
 B \frac{d^2 w}{d t^2} + \frac{d B}{d t} \frac{d w}{d t} = - D \frac{d w}{d t} .
\]
Differentiating once more and taking into account linearity of $B$ 
with respect to $t$ we obtain 
\[
 B \frac{d^3 w}{d t^3} + 2\frac{d B}{d t} \frac{d^2 w}{d t^2} = - D \frac{d^2 w}{d t^2} .
\]
Thus the  third order derivative of the solution can be written as
\[
  B \frac{d^3 w}{d t^3} = - \frac{2+\alpha }{\alpha } D \frac{d^2 w}{d t^2} .
\]
Substituting this relation into \eqref{28} we get the equation
\begin{equation}\label{29}
\begin{split}
 B_n \frac{w_{n+1} - w_{n-1}}{2\tau } & + D (\sigma w_{n+1} + (1-2\sigma) w_n + \sigma w_{n-1}) \\
 & = \left ( \sigma - \frac{2+\alpha }{6\alpha } \right ) \tau^2 D \frac{d^2 w}{d t^2}
 + \mathcal{O} (\tau^4) .
\end{split}
\end{equation}
We approximate the second order derivative in
\eqref{29} by the standard  central difference formula
\[
 \frac{d^2 w}{d t^2} = \frac{ w_{n+1} - 2w_n + w_{n-1}}{\tau^2}  +
 \mathcal{O} (\tau^2).
\]
Then from \eqref{29} we get the semi-difference scheme of the fourth 
approximation order  
\begin{equation}\label{30}
 B_n \frac{w_{n+1} - w_{n-1}}{2\tau } + D (\sigma_0 w_{n+1} + (1-2\sigma_0) w_n + \sigma_0 w_{n-1}) = 0,
\end{equation}
where the optimal weight  parameter $\sigma_0$ is given by
\begin{equation}\label{31}
 \sigma_0 = \frac{2+\alpha }{6\alpha } .
\end{equation}

\begin{thm}\label{t-3}
The three-level high order difference scheme  \eqref{30}, \eqref{31} 
is unconditionally stable with respect to the initial data.
\end{thm}
\begin{pf} The scheme \eqref{30}, \eqref{31} belongs to the family of three-level
weighted schemes \eqref{15}, \eqref{16}. Thus it is stable if $\sigma_0 > 0.25$.
This condition is satisfied for all  $0 < \alpha < 4$, including 
the considered fractional powers $0 < \alpha < 1$.  
\end{pf}

The implementation of the high order three-level difference scheme \eqref{30}, \eqref{31}
requires to specify the second  initial condition $w_1$. It should be computed 
with the same fourth order accuracy. We do not have any robust, unconditionally stable
and efficient two-level high-order difference scheme. 
In all computations presented in the next section
the initial condition $w_1$ is computed by using the symmetrical two-level 
scheme \eqref{12} and a sufficiently fine grid is constructed on the
time interval 
$[0, \tau]$. 

Let this interval be divided into $m$ sub-intervals. The approximate solutions $w_{\beta/m}$
are computed for time moments $t_{\beta/m} = \beta \tau /m$, $\beta = 1, \ldots, m$
by using the following scheme
\begin{equation}\label{32}
\begin{split}
 B_{(\beta+1/2)/m} & \frac{w_{(\beta+1)/m} - w_{\beta/m}}{\tau/m} + D \frac{w_{(\beta+1)/m} + w_{\beta/m}}{2}= 0,\\
 & \quad \beta = 0,1, \ldots, m-1.
\end{split}
\end{equation}
For sufficiently smooth solutions of the differential problem \eqref{8}  and
if $m \sim N$ then the solution $w_1$ is computed with the required accuracy 
$ \mathcal{O} (\tau^4)$.
We note that the computational complexity of the three-level algorithm is increased 
approximately twice if such  approach is applied to compute the initial condition $w_1$. 

\section{Numerical Experiments}

Here we present results of the numerical solution of a model problem
\eqref{1}, \eqref{2}, \eqref{3} in two spatial dimensions, where 
the computational domain is a unit square 
\[
 \Omega = \{\bm x:  \; \bm x= (x_1, x_2), \quad 0 < x_1  < 1, \quad  0 < x_2  < 1 \}.
\]
The coefficients of the operator
$\mathcal{A}$ and the right-hand side function $f$ in the equation
\eqref{4} are defined as
\[
 k(\bm x) = 1,
 \quad c(\bm x) = \left \{ \begin{array}{cc}
 100,   & x_1^2 + x_2^2 \leq 0.25 ,  \\[3pt]
 1,   &  x_1^2 + x_2^2 > 0.25 ,  
\end{array}
 \quad \mu(\bm x) = 0,
 \quad f(\bm x) = 1.
 \right.
\]
The piecewise linear continuous $P_1$ Lagrange elements are used to 
approximate the elliptic operator. The domain $\Omega$ is covered by 
the uniform grid with $50$ intervals in each direction. 

\begin{figure}[htp]
\begin{center}
    \includegraphics[width=0.65\linewidth] {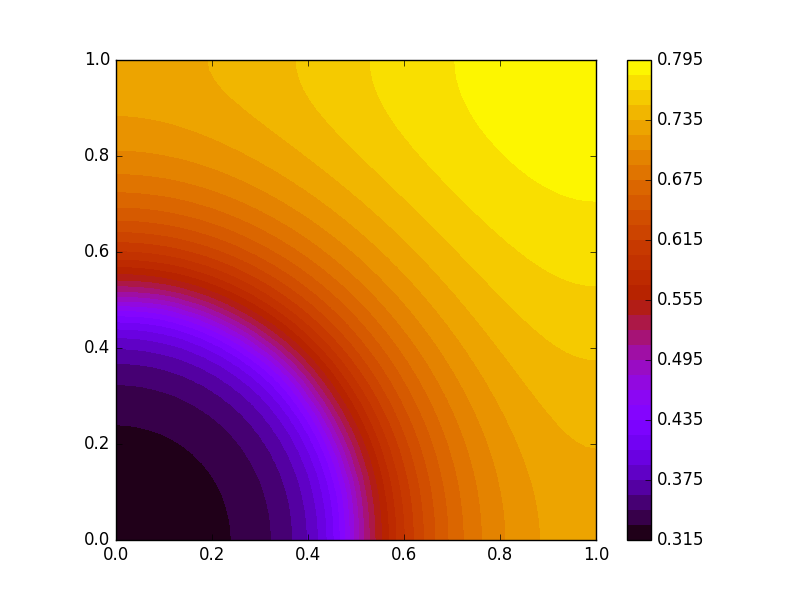} \\
    \includegraphics[width=0.65\linewidth] {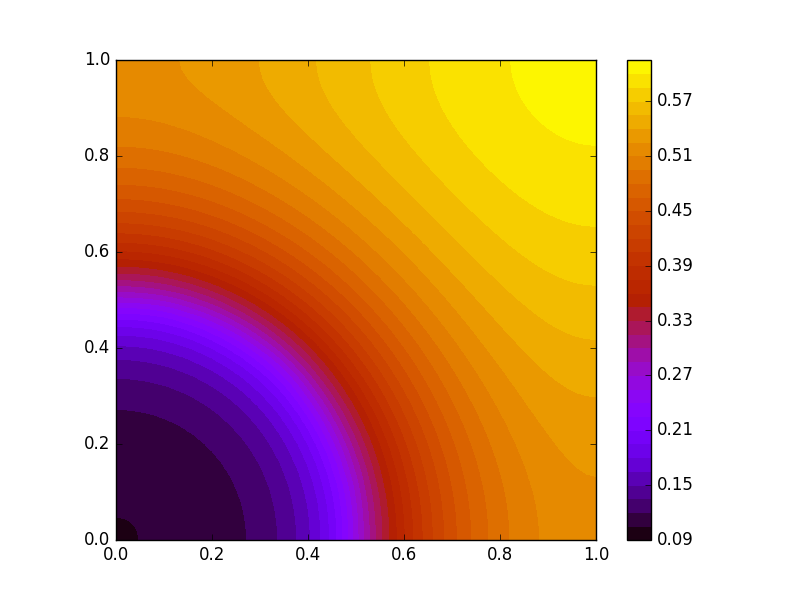} \\
    \includegraphics[width=0.65\linewidth] {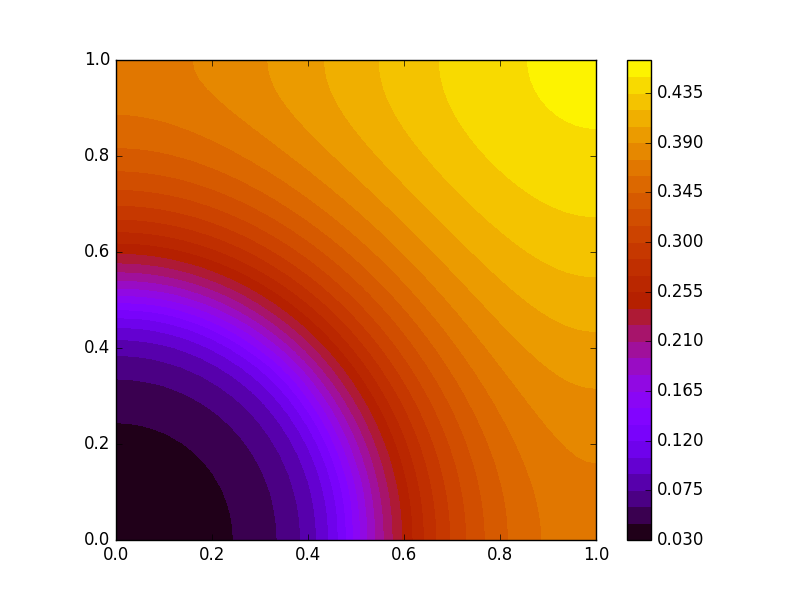} \\
\caption{The reference solution $\bar{w}_N$: top -- $\alpha = 
0.25$, center -- $\alpha = 0.5$,
         bottom -- $\alpha = 0.75$}
        \label{f-1}
\end{center}
\end{figure}

The accuracy of different approximations in time will be estimated by a reference solution.
It was obtained using the symmetrical two-level scheme  \eqref{12} with
$\sigma=0.5$ and taking  
a sufficiently small time step: $N = 5000$.
The relative errors of the approximate solution in the 
norm  $\|\cdot\|$ of space
 $L_2(\Omega)$ and in the norm  $\|\cdot\|_1$ 
$H^1(\Omega)$ are defined by
\[
 \varepsilon_1 = \frac{\|w_N - \bar{w}_N\|}{\|\bar{w}_N\|} ,
 \quad \varepsilon_2 = \frac{\|w_N - \bar{w}_N\|_1}{\|\bar{w}_N\|_1}, 
\]
where $\bar{w}_N$ is the reference solution.  
In Fig.~\ref{f-1}  we show the reference solution  for various values of the 
fractional power parameter $\alpha$. 

\begin{figure}[htp]
\begin{center}
    \includegraphics[width=0.8\linewidth] {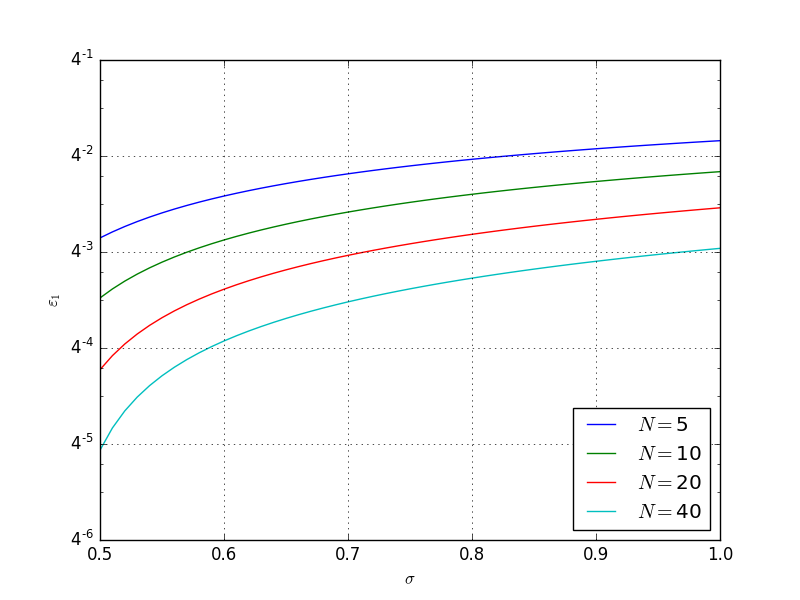} \\
    \includegraphics[width=0.8\linewidth] {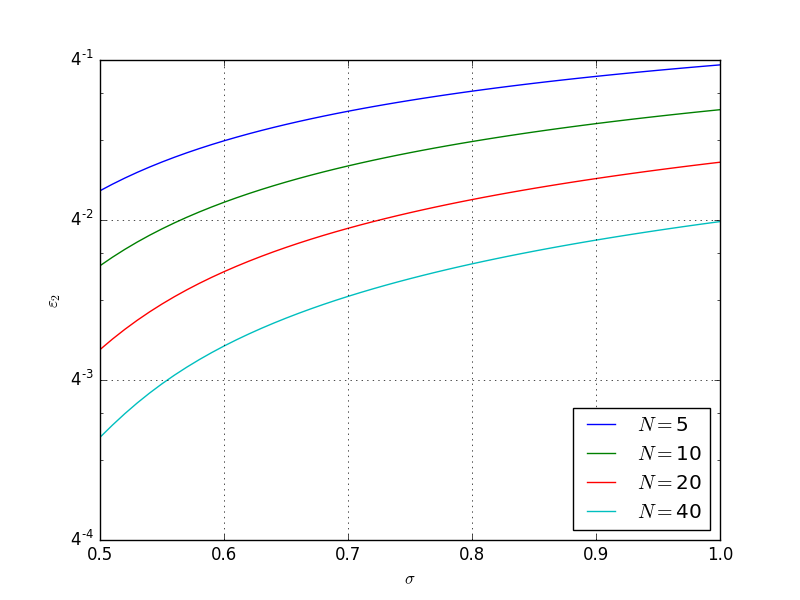} \\
\caption{The accuracy of the two-level scheme \eqref{12} for 
$\alpha = 0.25$: top -- $L_2(\Omega)$,
        bottom -- $H^1(\Omega)$}
        \label{f-2}
\end{center}
\end{figure}

\begin{figure}[htp]
\begin{center}
    \includegraphics[width=0.8\linewidth] {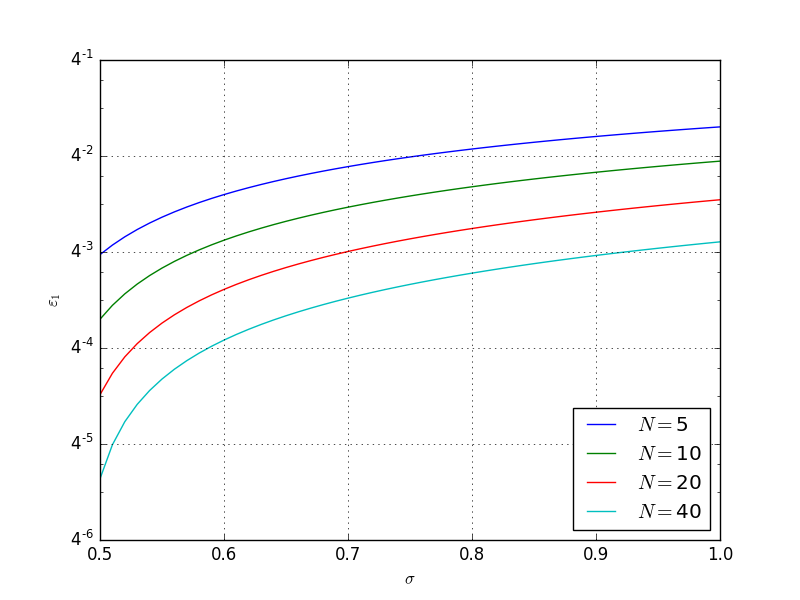} \\
    \includegraphics[width=0.8\linewidth] {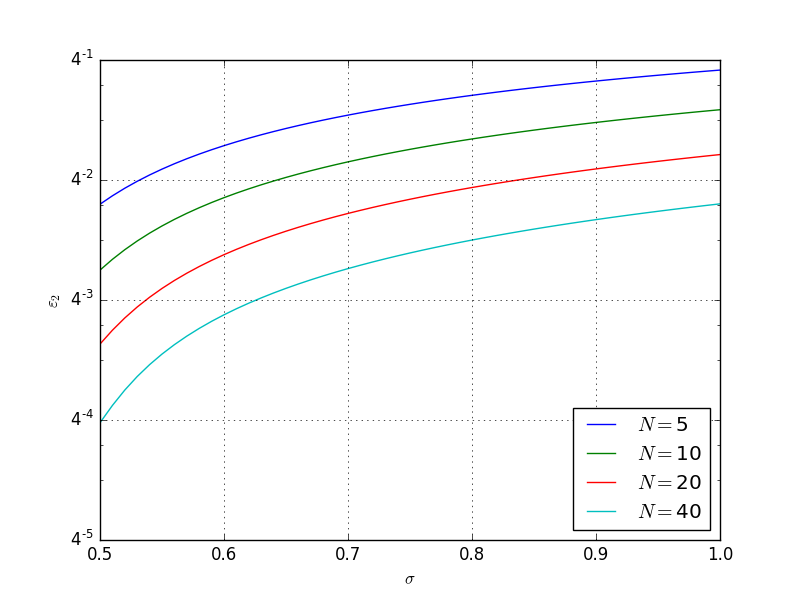} \\
\caption{The accuracy of the two-level scheme \eqref{12} for 
$\alpha = 0.5$:  top -- $L_2(\Omega)$,
        bottom -- $H^1(\Omega)$}
        \label{f-3}
\end{center}
\end{figure}

\begin{figure}[htp]
\begin{center}
    \includegraphics[width=0.8\linewidth] {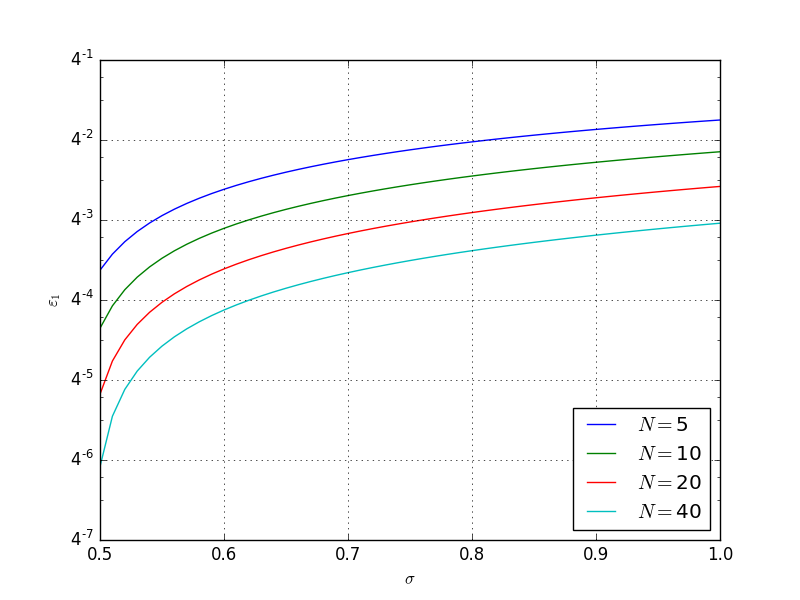} \\
    \includegraphics[width=0.8\linewidth] {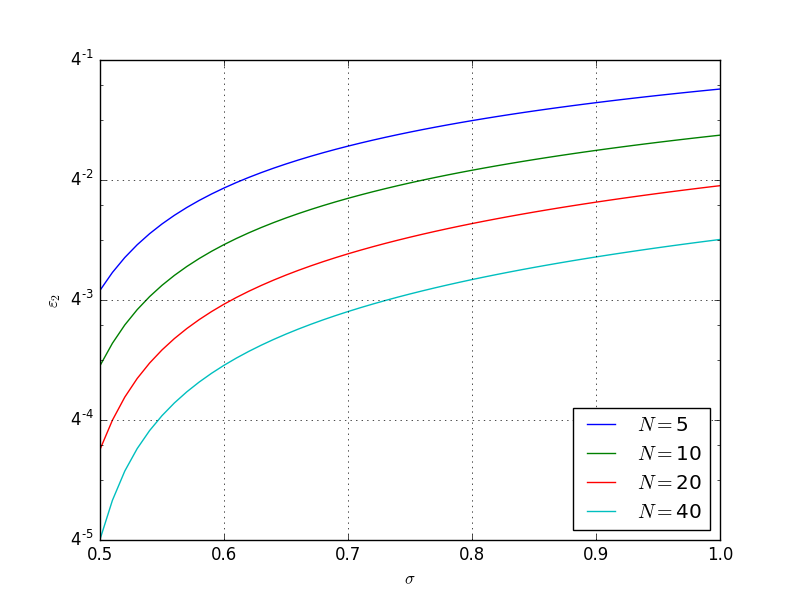} \\
\caption{The accuracy of the two-level scheme \eqref{12} for 
$\alpha = 0.75$:  top -- $L_2(\Omega)$,
        bottom -- $H^1(\Omega)$}
        \label{f-4}
\end{center}
\end{figure}

For the two-level weighted difference scheme \eqref{12} the  errors of the
solution are presented in Figs.~\ref{f-2} -- \ref{f-4}. As it 
follows from the theoretical analysis the accuracy of approximation
is essentially increased 
for the values of parameter $\sigma$ in the
neighbourhood of 0.5. 

We also note that the accuracy of the approximate
solution is better for larger values of $\alpha$. This result  is explained 
by the increased smoothness of the solution for larger values of $\alpha$.

\begin{figure}[htp]
\begin{center}
    \includegraphics[width=0.8\linewidth] {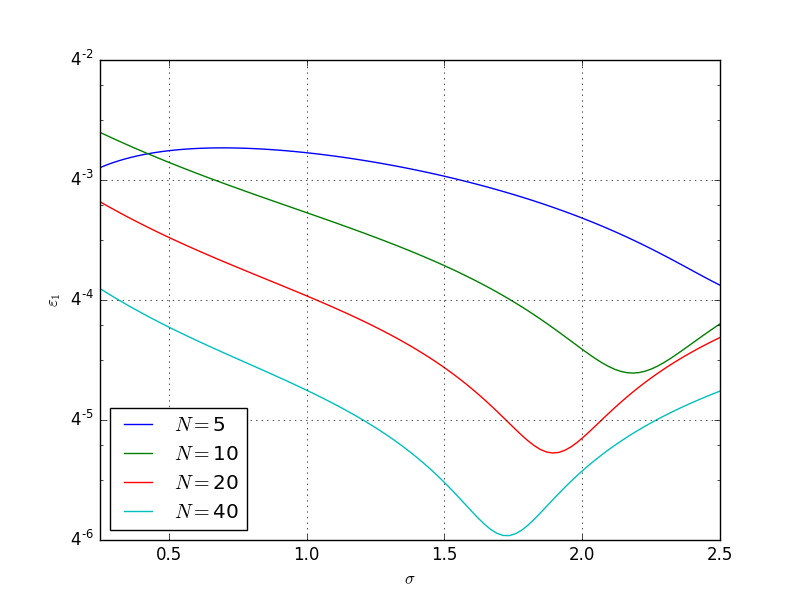} \\
    \includegraphics[width=0.8\linewidth] {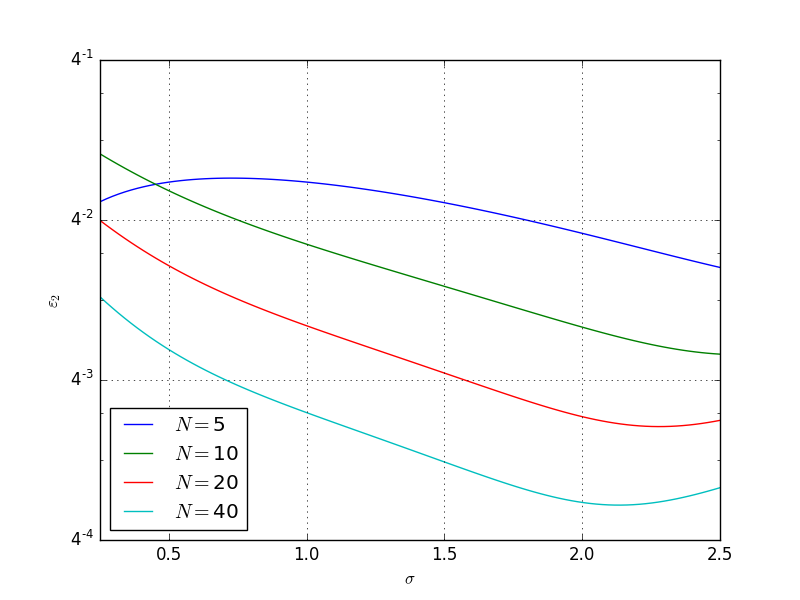} \\
\caption{The accuracy of the three-level scheme \eqref{15} with 
the initial condition computed using the standard two-level scheme  
\eqref{21} for $\alpha = 0.25$ ($\sigma_0 = 3/2$): 
top -- $L_2(\Omega)$,
        bottom -- $H^1(\Omega)$.}
        \label{f-5}
\end{center}
\end{figure}

\begin{figure}[htp]
\begin{center}
    \includegraphics[width=0.8\linewidth] {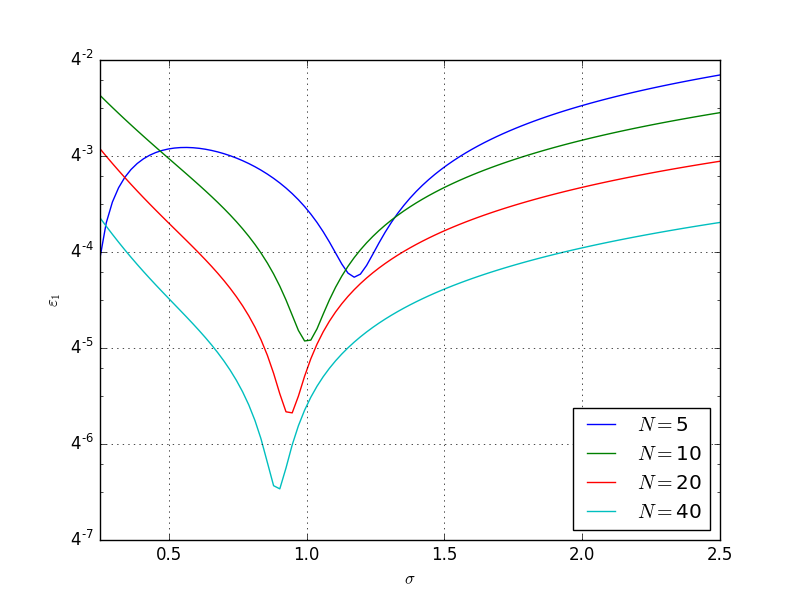} \\
    \includegraphics[width=0.8\linewidth] {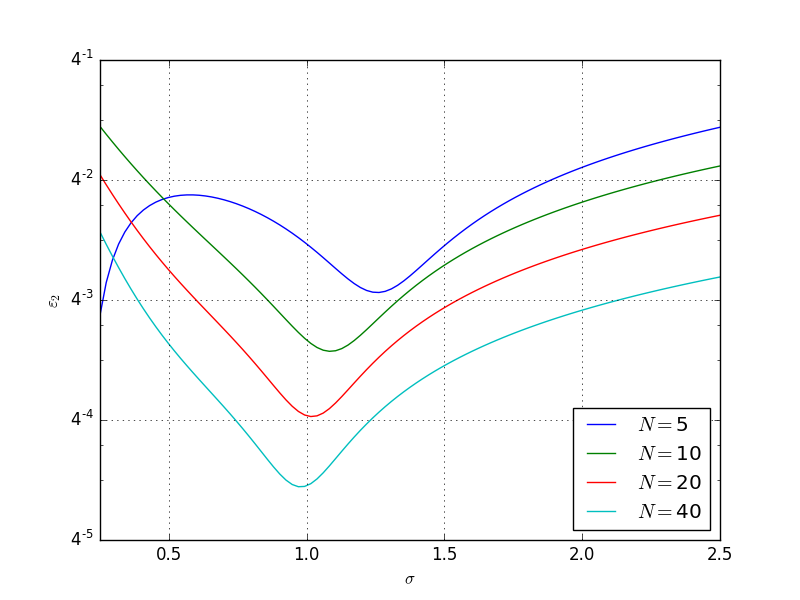} \\
\caption{The accuracy of the three-level scheme \eqref{15} with 
the initial condition computed using the standard two-level scheme  
\eqref{21} for $\alpha = 0.5$ ($\sigma_0 = 5/6$): 
top -- $L_2(\Omega)$,
        bottom -- $H^1(\Omega)$.}
        \label{f-6}
\end{center}
\end{figure}

\begin{figure}[htp]
\begin{center}
    \includegraphics[width=0.8\linewidth] {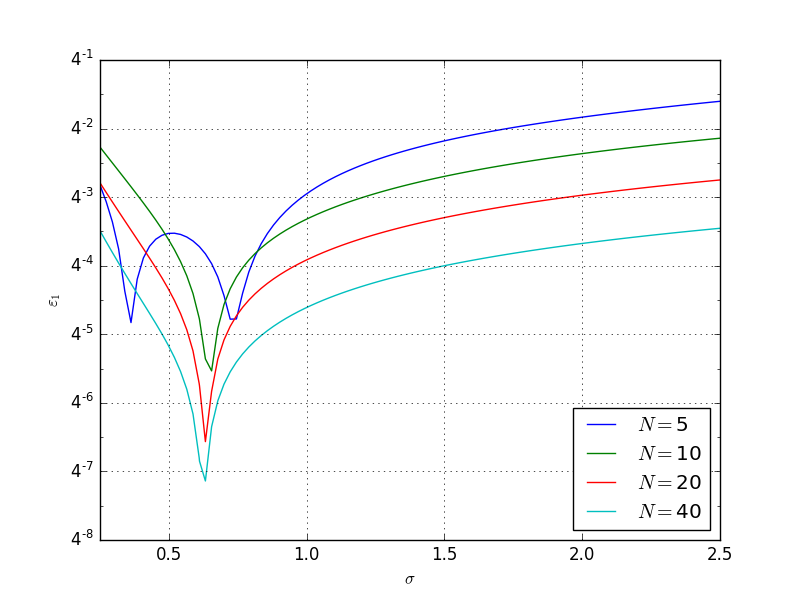} \\
    \includegraphics[width=0.8\linewidth] {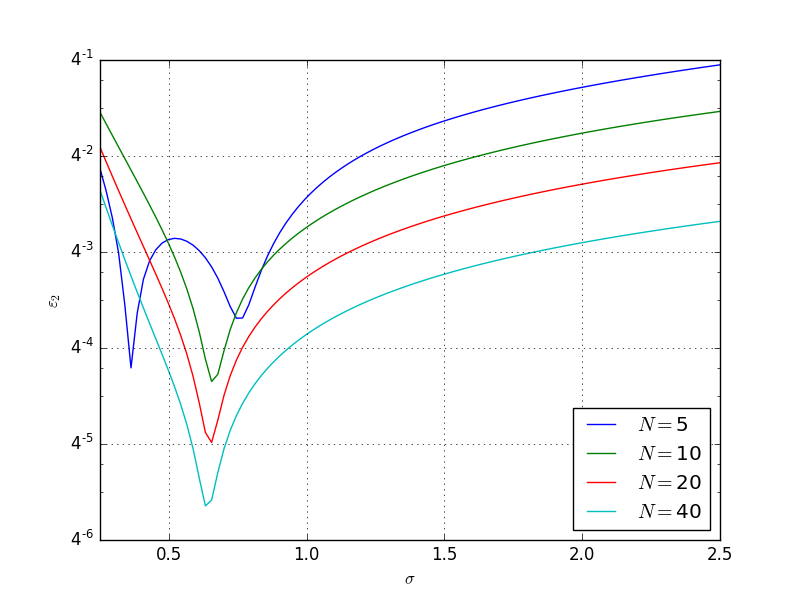} \\
\caption{The accuracy of the three-level scheme \eqref{15} with 
the initial condition computed using the standard two-level scheme  
\eqref{21} for $\alpha = 0.75$ ($\sigma_0 = 11/18$): 
top -- $L_2(\Omega)$,
        bottom -- $H^1(\Omega)$.}
        \label{f-7}
\end{center}
\end{figure}

The  main goal of this paper is to investigate the accuracy of the 
three-level difference scheme \eqref{15}, \eqref{16}.
  First we have
used the two-level symmetrical scheme \eqref{21} to compute the
initial condition for $w_1$. For sufficiently smooth solutions 
it defines the initial condition with $O(\tau^2)$ accuracy. Results 
of computations for various weight parameters $\sigma$ are shown in      
Figs.~\ref{f-5} -- \ref{f-7}. It is clearly seen that the accuracy 
of the approximation is increased for the optimal weight parameter $\sigma_0$.

\begin{figure}[htp]
\begin{center}
    \includegraphics[width=0.8\linewidth] {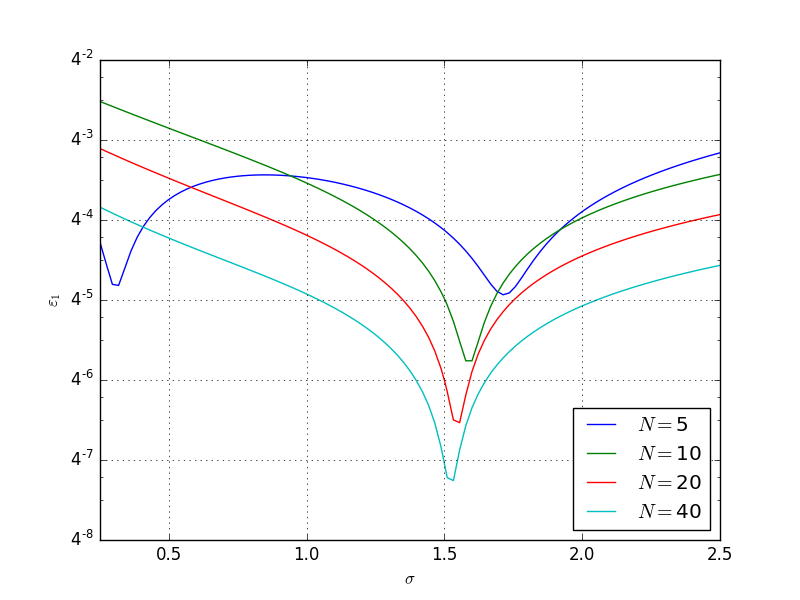} \\
    \includegraphics[width=0.8\linewidth] {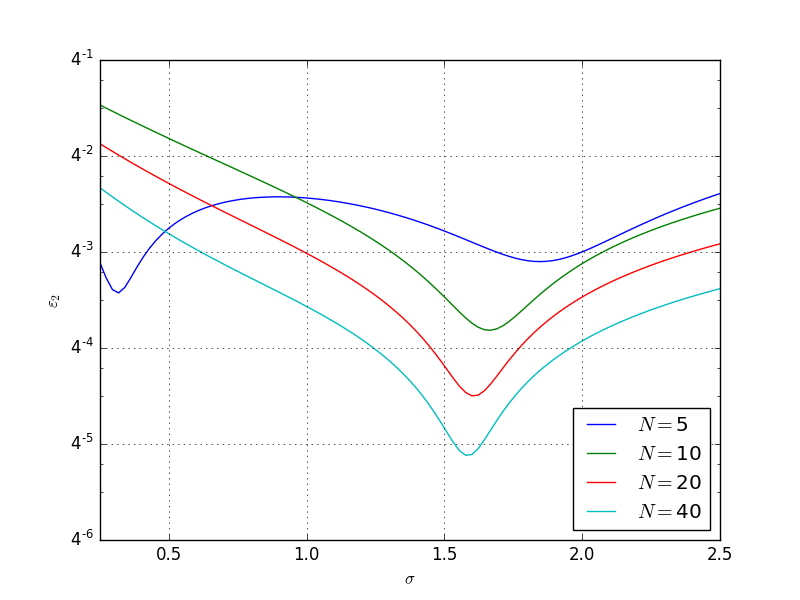} \\
\caption{The accuracy of the three-level scheme \eqref{15} with 
the initial condition computed using the two-level scheme  
\eqref{32} for $\alpha = 0.25$ ($\sigma_0 = 3/2$): 
top -- $L_2(\Omega)$,
        bottom -- $H^1(\Omega)$.}
        \label{f-8}
\end{center}
\end{figure}

\begin{figure}[htp]
\begin{center}
    \includegraphics[width=0.8\linewidth] {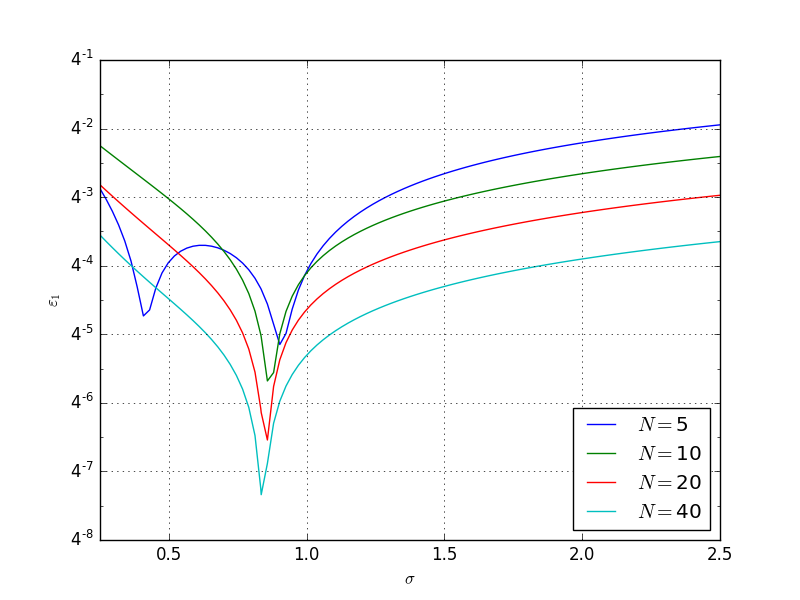} \\
    \includegraphics[width=0.8\linewidth] {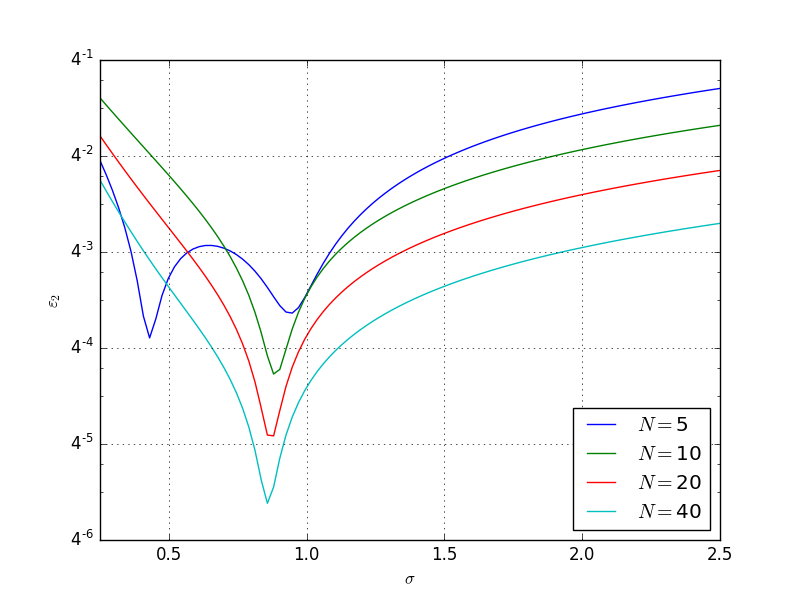} \\
\caption{The accuracy of the three-level scheme \eqref{15} with 
the initial condition computed using the two-level scheme  
\eqref{32} for $\alpha = 0.5$ ($\sigma_0 = 5/6$): 
top -- $L_2(\Omega)$,
        bottom -- $H^1(\Omega)$.}
        \label{f-9}
\end{center}
\end{figure}

\begin{figure}[htp]
\begin{center}
    \includegraphics[width=0.8\linewidth] {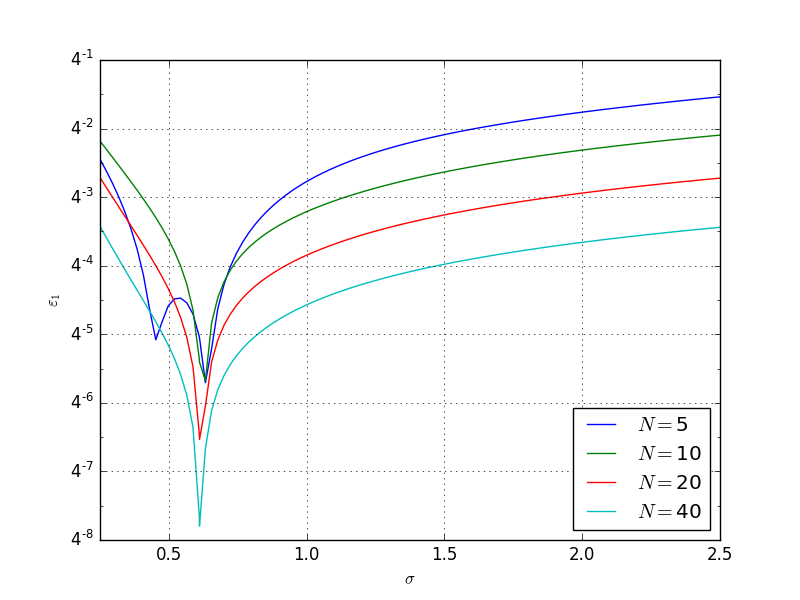} \\
    \includegraphics[width=0.8\linewidth] {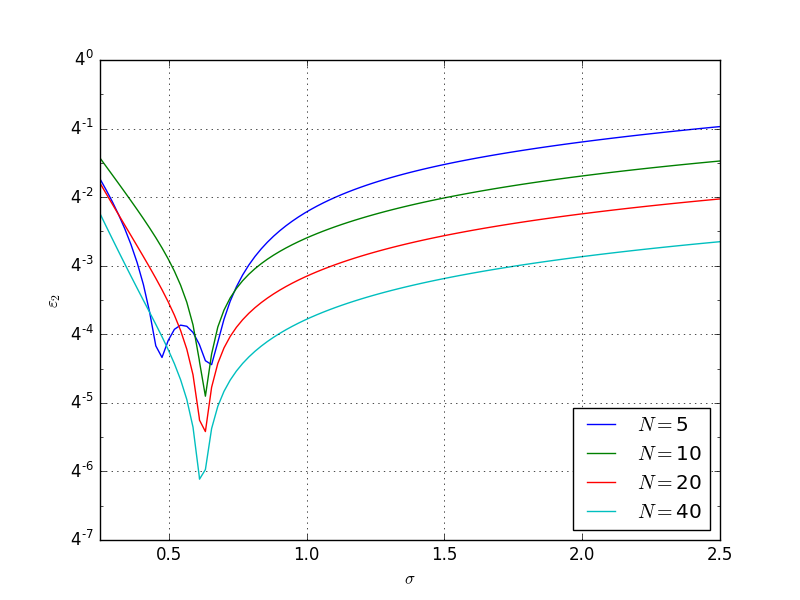} \\
\caption{The accuracy of the three-level scheme \eqref{15} with 
the initial condition computed using the two-level scheme  
\eqref{32} for $\alpha = 0.75$ ($\sigma_0 = 11/18$): 
top -- $L_2(\Omega)$,
        bottom -- $H^1(\Omega)$.}
        \label{f-10}
\end{center}
\end{figure}

Next we have investigated the influence of the initial condition for $w_1$. 
The accuracy of the approximation  is further increased when the initial condition 
is computed using the algorithm \eqref{32}. Results of computations for various 
weight parameters $\sigma$ are shown in
Figs.~\ref{f-8} -- \ref{f-10}. 

Here we note that 
the observed convergence rates of the two-level and three-level  schemes  
depend on the discrete regularity of the solution of the discrete 
fractional power problem and they are not reaching the maximal possible 
convergence rates of these schemes. 
As expected from the theoretical analysis (see, e.g.    
\cite{duan2018numerical}), the convergence rate is increased for larger 
values of $\alpha$. The  dependence on the regularity of the
solution  can be reduced by using 
geometrically refined time grids. 

\section{Conclusions}

1. We have formulated the problem of finding the high order difference schemes for solving 
the nonstationary Cauchy type problem which is equivalent to the fractional power elliptic 
problem. The high order approximations are used to approximate the time dependence 
of the solution, while the elliptic operator is approximated by the standard finite
element scheme. 

2. The sufficient stability conditions are given for the two-level discrete schemes
with weight parameters. The second order accuracy is proved for the symmetrical 
Crank-Nicolson type scheme. 

3. The family of three-level symmetrical discrete schemes is
constructed and investigated. It is proved 
that the second order approximation is valid for sufficiently smooth solutions. The initial 
condition on the first time level is computed by using the symmetrical two-level
scheme.   

4. It is shown that for a special weight parameter $\sigma_0$ we get the fourth-order 
three-level scheme. The value of this optimal parameter depends on the fractional
power $\alpha$ of the elliptic operator. 
The initial 
condition on the first time level of the main grid is computed by using the 
symmetrical two-level
scheme with a specially selected fine time grid. 

5. The theoretical results are illustrated by results of numerical experiments. 
A two-dimensional problem is solved for the elliptic operator with  the discontinuous 
sink term  coefficient.

\section*{Acknowledgements}
This work of second author was supported by the mega-grant of the Russian Federation Government (\#~14.Y26.31.0013).

\end{document}